\documentclass[12pt]{amsart}

\usepackage{amsmath,amscd,amssymb,amsfonts}
\usepackage{xcolor}
\usepackage{graphicx}
\usepackage{pb-diagram}

\newcommand{\Z}{\mathbb{Z}}
\newcommand{\Q}{\mathbb{Q}}

\newcommand{\C}{\mathbb{C}}

\renewcommand{\H}{\mathbb{H}}
\renewcommand{\P}{\mathbb{P}}

\newcommand{\XS}{\mathcal{X}}
\newcommand{\YS}{\mathcal{Y}}
\newcommand{\ZS}{\mathcal{Z}}

\newcommand{\lra}{\longrightarrow}

\title{A special Calabi-Yau degeneration with trivial monodromy}
\author{Slawomir Cynk and Duco van Straten}
\begin{document}
\begin{abstract}
A well-known theorem of Kulikov, Persson and Pinkham states that a 
degeneration of a family of K3-surfaces with trivial monodromy can
be completed to a smooth family. We give a simple example that an
analogous statement does not hold for Calabi-Yau threefolds.
\end{abstract}
\maketitle

\section{Introduction}

In the study of degenerations of $K3$ surfaces, the theorems
of Kulikov \cite{Kulikov} Persson and Pinkham \cite{PerssonPinkham} 
are fundamental and play a key role in the proof of the surjectivity 
of the period map for $K3$ surfaces. One important result concerns
{\em degenerations of type I}: if $f: \XS \lra \Delta$ is a degeneration 
of $K3$ surfaces over the disc $\Delta:=\{ t \in \C\;|\;|t| <1\}$, 
with monodromy  on $H^2(X_t)$ of {\em finite order}, then after an 
appropriate base change and birational modification of the zero-fibre 
we obtain a family  for which the fibre $X_0$ over $0$ is {\em smooth}. 
For example, if $\XS \lra \Delta$ acquires only ADE-singularities at the
zero-fibre, then the theorem is implied by to the phenomenon 
of {\em simultaneous resolution after base change} for these singularities, 
discovered by Brieskorn \cite{Brieskorn}, generalising the $A_1$-case 
of Atiyah \cite{Atiyah}. The fact that for ADE-singularities in dimension 
two resolving and deforming 'are the same' has been recognized as a typical 
feature in hyperk\"ahler 
geometry and recently a generalisation of the above theorem to 
degenerations of higher dimensional  hyperk\"ahler manifolds was  
given in \cite{KLSV}. As remarked in that paper, it is known that this does 
not generalise to degenerations of higher dimensional {\em Calabi-Yau} 
varieties. The homological monodromy of an odd dimensional $A_1$-singularity 
has infinite order, but the $A_2$-singularity has finite order. It was 
remarked long ago by Clemens, Friedman \cite{Friedman2} and Morgan 
\cite{Morgan} that a Calabi-Yau $3$-fold degeneration acquiring an $A_2$-singularity does not have a smooth filling and there is a result of Voisin 
\cite{Voisin} implying a same result for four dimensional varieties aquiring 
an $A_1$-singularity.\\

The purpose of this paper is to present a completely different type of example
and analyse it in some detail. The main result is the following\\

\centerline{\Large \bf \em Theorem} 
\vskip 20pt
{\em There exists a flat family $f:\YS \lra \Delta$ of
projective threefolds such that
\begin{enumerate}
\item $\mathcal Y$ is a smooth fourfold,
\item for $t \neq 0$, the fibre $Y_{t}:=f^{-1}(t)$ is a smooth 
Calabi-Yau threefold with
  $$h^{1,1}(Y_{t})=41,\;\;\;\; h^{1,2}(Y_{t})=1,$$
\item the singular locus of $Y_{0}$ is a line $L$, $Y_0$ is double along
  $L$ with exactly four pinch points, 
\item the elliptic curve $E$ doubly covering $L$ and ramified over the set
$\Sigma$ of these four pinch points has $j$-invariant equal to $1728$,
\item the blow-up of $Y_{0}$ in $L$ is a smooth Calabi-Yau threefold $Z_0$
  with $$h^{1,1}=46,\;\;\;h^{1,2}=0,$$
\item the local system $H^{i}(Y_{t})$ has trivial monodromy over
    $\Delta^{*}$ for $i\not=3$ and $\Z/2\Z$-monodromy for $i=3$.
\end{enumerate}
}
\centerline{\bf \em Corollary} 
\vskip 10pt
{\em The semi-stable reduction, obtained after a base change 
$t \mapsto t^2$,
has trivial monodromy and the special fibre consists of two components,
one of which is a projective, smooth and rigid Calabi-Yau manifold and the
other is a smooth quadric bundle.}\\

In the Clemens-Friedman-Morgan example we have a Calabi-Yau variety that 
acquires an $A_2$-singularity, which is a {\em terminal singularity}. So 
from the point of view of the minimal model program, this singular variety 
should be considered as good as a smooth one. In fact, there does not exist 
a crepant resolution of the singular  member. 
In our example the zero-fibre $Y_0$ has {\em canonical singularities}, 
and admits a crepant resolution to a honest smooth Calabi-Yau variety $Z_0$. 
There is a change in cohomology in going from $Z_0$ to $Y_t$. 
In section $4$ we will see that the third limiting 
Hodge structure splits
\[ H^3_{lim}(Y_{\infty},\Q) = H^3(Z_0,\Q) \oplus H^1(E,\Q)(-1)\]
where $E$ is the elliptic curve mentioned in the theorem. The two pieces 
of this decomposition are necessarily supported by two different irreducible 
components of the zero fibre.\\

The structure of the paper is as follows. In section 1 we describe the
basic structure of the example. In section 2 we describe the resolution 
process in some detail and give the proof of the above theorem. 
In section 3 we analyse the example on a cohomological level. 
In section 4 we collect some remarks and speculations.\\

\section{The example}
\subsection{Double octics}
Our example is based on certain special {\em double octics}, 
double covers of $\P^3$ ramified along an arrangement of eight planes 
$P_1,P_2,\ldots, P_8$. Such a space can be given as a hypersurface 
in a weighted projective space 
\[ \{u^2=L_1 L_2 L_3 L_4 L_5 L_6 L_7 L_8 \} \subset \P[4,1,1,1,1], \]
where $L_i$ is a linear form defining the plane $P_i$. 
For a generic choice of the planes $P_i$, the branch divisor $D=\cup_{i=1}^8 P_i$
(and hence also the double cover) has $28$ double lines 
and $56$  triple points, along which it is singular. By blowing up the 
(strict transforms of) the double lines in any order we obtain a crepant 
resolution which is a Calabi -- Yau threefold with  Hodge numbers
$h^{12}=9, h^{11}=29$.
For special positions of the planes the singularities of the double octic
change, but as long as the configuration of planes does not have $4$-fold
lines or $6$-fold points, there still exists a crepant resolution, but now 
the Hodge numbers can take various values, depending on the precise properties 
of the configuration. Recently, all different cases have been listed.
For more information on special double octics we refer to 
\cite{Meyer}, \cite{CynkKocel}.\\ 
 
Our example is based on the following family of double octics
\[X_t:=\{u^{2}=xy \left( x+y \right) z \left( x+2\,y+z +tv\right)
 v \left( y+z+v \right) \left( x+y+z+ \left( t-1 \right) v \right)\}. \]
where $t \in \P^1$ is considered as the parameter. 
We label the linear forms $L_1,L_2,\ldots,L_8$ in the order as they are 
written in the above equation:
\[ L_1=x,\;\; L_2=y,\;\;L_3=x+y,\ldots,\;\; L_8=x+y+z+\left( t-1 \right) v .\]

For $t \neq 0,1,2,\infty$ this configuration has exactly

\begin{itemize}
\item a single triple line $\ell_{triple}: x=y=0$,
\item 25 double lines,
\item six fourfold points, 
\[(1:0:0:0),\;\;(0:1:-1:0),\;\;(1:-1:1:0),\;\;(t-2:1:0:-1),\;\;(1:0:-1:0),\;\;(1:-1:0:0)\]
not lying on the triple line $\ell_{triple}$ 
(called points of \emph{type} $p_{4}^{0}$),
\item five fourfold points 
\[ (0: 0: 1: 0),\;\;(0: 0: 0: 1),\;\;(0: 0: -t: 1),\;\;(0: 0: 1: -1),\;\;
(0: 0: t-1: -1)\] 
on the triple line $\ell_{triple}$ (called points of \emph{type} $p_{4}^{1}$).
\end{itemize}

This arrangement is projectively equivalent to arrangement No. $153$
in \cite{Meyer} via the coordinate transformation
\[(x,y,z,v)\mapsto(-y-z-v,v,y,x+y+z+ \left( t+1 \right) v),\]
and reparametrisation $t \mapsto t-2$.

We are here concerned specifically with the degeneration that occurs at 
$t = 0$ and consider the family of double octics over the unit disc:
\[ \pi: \XS \to \Delta, \]
with fibre over $t \in \Delta$ the double octic $X_t$ defined by the 
above equation. For this degeneration, there are two important lines,
namely the triple line $\ell_{triple}$ and the moving line $m_t=P_4 \cap P_5$. 
If $t$ goes to $0$, the moving line intersects the triple line, and the two 
$4$-fold points $ \ell_{triple} \cap P_4$ and $\ell_{triple} \cap P_5$ move 
together to form a $5$-fold point. The result is that for $t=0$ we
obtain a configuration equivalent to the rigid arrangement No. $93$
of \cite{Meyer}. \\
\begin{center}
\includegraphics[height=6cm]{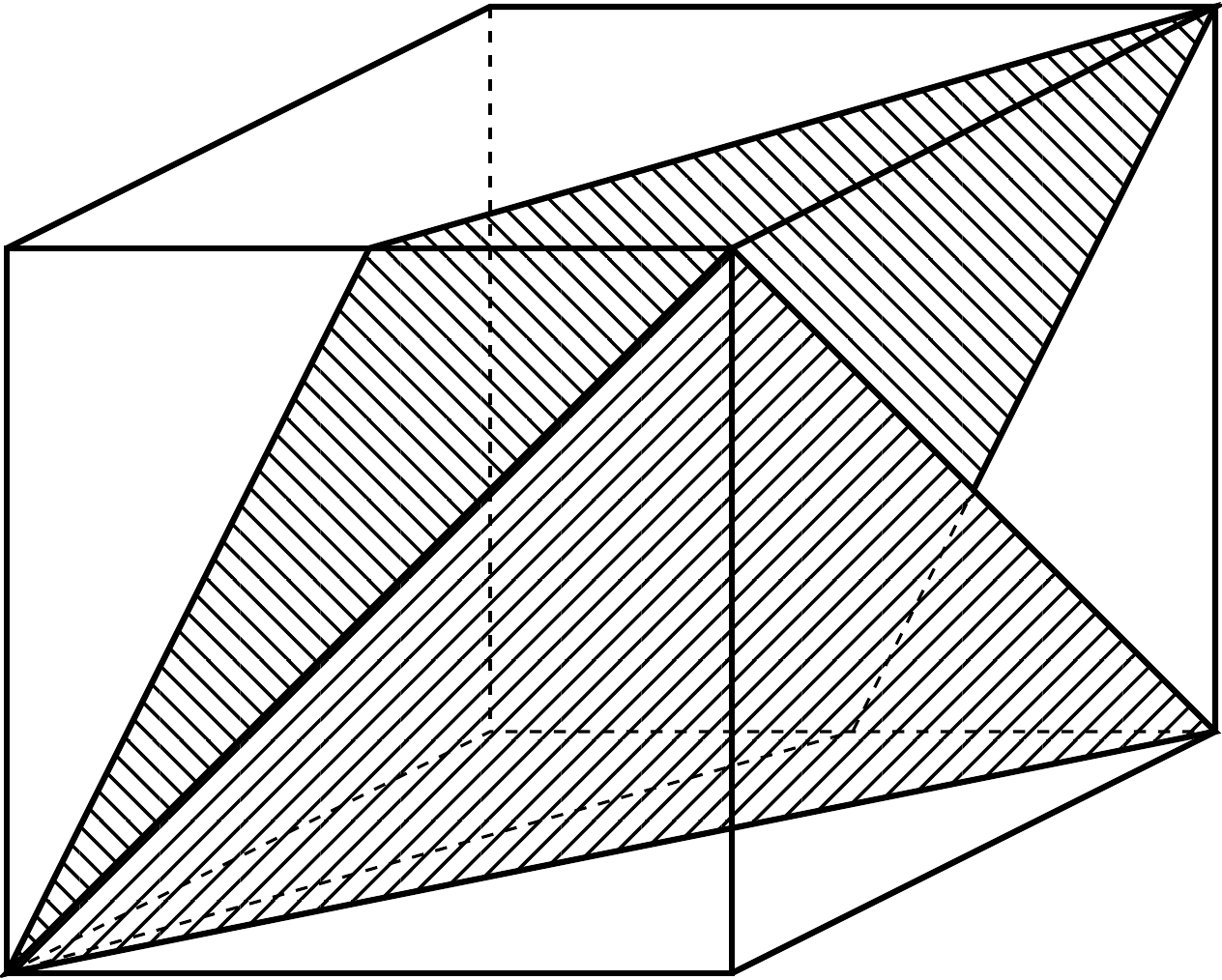}
\end{center}

\bigskip

The plane spanned by the lines $\ell_{triple}$ and $m_0$ makes a fourth 
plane $P$ through $\ell_{triple}$. This plane does not belong to
the octic arrangement, but the the planes $P_1,P_2,P_3, P$ define four 
points on the projective line $L$ of all planes through $\ell_{triple}$. 
We let $E$ be the double cover of $L$ ramified over these four points. 
As these planes are in harmonic position, the $j$-invariant of $E$
is seen to be $1728$. This is the elliptic curve we were alluding to in 
the introduction. We note furthermore that all strata of the singular
locus of $X_t$ are defined by the intersections among the planes $P_i$.
As these are not interchanged by the monodromy, these loci form trivial 
families over $\Delta^{*}$.
 
\subsection{Construction of $\YS \lra \Delta$}

We will now construct the family $f:\YS \lra \Delta$ of the theorem
from the family of singular double octics $\XS \lra \Delta$ by a 
certain specific sequence of blow-ups. To be more precise, we construct 
a sequence of blow-ups and a diagram of two-fold covers over it:

\[
  \begin{diagram}\dgARROWLENGTH=1cm
    \node{\XS^{(3)}}\arrow{e} \node{\XS^{(2)}}\arrow{e}\arrow{s}
    \node{\XS^{(1)}}\arrow{e}\arrow{s}\node{\XS}\arrow{s}\\ 
    \node[2]{\P^{(2)}}\arrow{e}\node{\P^{(1)}}\arrow{e}\node{\P^{3}\times\Delta}
  \end{diagram}
\]

We are dealing here with families over $\Delta$ and by {\em blowing-up 
a line} we mean the blowing-up a {\em relative} line over $\Delta$,
which is really a surface. 

\begin{itemize}
\item Blow-up in $\P^3 \times \Delta$ first in the locus of all fourfold points of type $p_{4}^{0}$ and the triple line. Call the resulting space $\P^{(1)}$ and let $D^{(1)}$ be the strict transform of $D=D^{(0)}$, plus the 
divisor $E:=L \times \ell_{triple} \times \Delta = \P^1 \times \P^{1} \times \Delta$ 
lying over the triple line $\ell_{triple} \times \Delta$ over $\Delta$. 
We denote by $\XS^{(1)}$ the double cover ramified over $D^{(1)}$. 
The branch divisor $D^{(1)}$ contains, apart from the strict transforms of the
$25$ double lines of $D$ also  eight 'new' double lines: 
three lines $m_{1}, m_{2}, m_{3}$ in one ruling corresponding to the 
three planes $P_1,P_2,P_3$ containing the triple line $\ell_{triple}$, 
and five further lines $m_{4},\dots, m_{8}$ in the second ruling, namely 
the intersection of the strict transforms of the planes $P_4,P_5,\ldots, P_8$.
In local coordinates near the line $m_1$ the space $\XS^{(1)}$ is described
by the equation
\[u^{2}=xyz(x+2xy+z+t)F.\]
Here the divisor $x=0$ is the exceptional set $E$, $y=0$ defines the line
$m_1$. The next two factors $z=0$ and $x+2xy+z+t=0$ are the equations
of the strict transforms of  $L_4$ resp $L_5$. The factor $F$ is the product
of all other factors; it is non-zero on $m_1$.
(We shall consider only the resolution of the complement of lines $m_{2}$ and $m_{3}$, as the resolution in neighbourhoods of these two lines is completely
analogous to what happens near $m_1$).

\item In the next step we blow-up all double curves, except three of them: 
$m_{4},m_{5}$ and the intersection of strict transforms of 
planes $P_{4}$ and $P_{5}$. The new branch  divisor $D^{(2)}$ is the strict 
transform of $D^{(1)}$. Call the resulting double cover $\XS^{(2)}$.
Over the blow - up of the line $m_1$ the space $\XS^{(2)}$ is
described in two charts by
\vskip 10pt
\begin{itemize}
\item [] $(x,y)\longmapsto (xy,y)$,\qquad $u^{2}=xz(xy+2xy^{2}+z+t)F.$
\item [] $(x,y)\longmapsto (x,xy)$,\qquad $u^{2}=yz(x+2x^{2}y+z+t)F$
\end{itemize}
\vskip 10pt
\item The last step of the resolution is to blow-up 
the double cover in the union of singular loci of all fibers
$\bigcup_{t\in\Delta}\operatorname{Sing}(X_{t})$. 
We denote by $\YS=\XS^{(3)}$ the resulting variety. To analyse what we
have, we look in the two charts of $\XS^{(2)}$ around the blow-up 
of $m_1$.  In the second chart the branch-divisor and the fibre 
of the branch-divisor are both simple normal crossings. As a
consequence, the  parts of the spaces $\YS$ and $Y_{t}$ lying over this
chart are smooth. To analyse $\YS$ in the first chart, we have to blow-up
$\XS^{(2)}$ in the ideal
\[(xz, x(xy+2xy^{2}+z+t), z(xy+2xy^{2}+z+t),u). \]
This blow-up is given as the closure of the map
\[(x,y,z,t,u) \mapsto (X,Y,Z,T)=(xz,x(xy+2xy^{2}+z+t), z(xy+2xy^{2}+z+t),u)\]
Using {\sc Singular} \cite{Singular} we verified that $\YS$ is smooth in this
chart as well and the special fibre is singular along the line
\[\operatorname{Sing}(Y_{0})=(x=z=u=X=Y=Z=0)\]
is contained in the affine chart $T=1$, moreover in this affine chart
the variety $Z$ is locally given as
\[x=XY, z=XZ, u=XYZ, (y+2y^{2})XY+XZ-YZ+t=0.\]
\end{itemize}

{\bf proof of the Theorem:} 

{\em Hodge numbers:} The statement about the Hodge numbers follow from the formulas that express the Hodge-numbers of a resolved double octics $X_t$ and are recorded in \cite{Meyer}: $h^{1,1}=46, h^{1,2}=0$ for $X_0$ and $h^{1,1}=41$, $h^{1,2}=1$ for $X_t$, $t \in \Delta^*$. As the Hodge numbers do not depend on the
choice of the resolution, we find the Hodge numbers as stated in the theorem.\\

{\em Properties of $Y_0$:} The fact that $Y_0$ is double along the
line $L$ and is resolved by a single blow--up follows from the above 
local calculations. The singular line $L$ of $Y_{0}$ can be identified with
the pencil of planes through the triple line $\ell_{triple}$ at the fivefold 
point $(0:0:0:1)$. Three of the pinch-points correspond to three planes 
containing $\ell_{triple}$, the fourth pinch-point is the direction of the intersection line $P_{4} \cap P_{5}$.\\ 

{\em Monodromy:} The cohomology group  $H^{2}(Y_{t},\C), t\in\Delta^{*}$ is 
generated by classes of components of the exceptional locus of the crepant
resolution, so it has trivial monodromy. 
A simple way to determine the monodromy on $H^3(Y_t)$ is using the
{\em Picard-Fuchs operator $\mathcal{P}$} of the family that was 
determined in \cite{CvS}. 
\begin{eqnarray*}
&&\mathcal{P}=4\Theta(\Theta-\tfrac12)(\Theta-\tfrac32)(\Theta-2)
  -12t(\Theta-\tfrac12)^{2}(\Theta^{2}-\Theta+\tfrac1{12})
  +13t^{2}(\Theta^{4}+\tfrac9{26}\Theta^{2}+\tfrac1{208})\\
&&  -6t^{3}(\Theta+\tfrac12)^{2}(\Theta^{2}+\Theta+\tfrac7{13})
  +t^{4}(\Theta+\tfrac12)(\Theta+1)^{2}(\Theta+\tfrac32) \in \Q \langle t,\Theta \rangle,  
\end{eqnarray*}
where $\Theta=t\partial/\partial t$ is the logarithmic derivative with respect to the parameter $t$. The Riemann symbol, that collects the exponents at all
singular points is:
\[
\left\{ 
\begin {array}{cccc}
2&1&0&\infty\\
\hline
0&0&0&1/2\\
1/2&1/2&1/2&1\\
1/2&1/2&3/2&1\\
1&1&2&3/2\\
\end {array} \right\} 
\]
The local system $\mathcal{S}ol$ of solutions on $\P^1 \setminus \{0,1,2,\infty \}$ is isomorphic to the local system with fibres $H^3(Y_t,\C)$. 
It can be checked by computing the formal solutions around $0$ that
{\em no logarithms occur in the solutions} and one finds four series solutions
\[\phi_0(t)=1+\tfrac14t-\tfrac{23}{1120}t^{2}+\ldots,\;\;\; \phi_1(t)
  =t^{1/2}(1+t+\tfrac 45t^{2}\ldots), \]
\[ \phi_2(t)=t^{3/2}(1+t+\tfrac67t^{2}+\ldots),\;\;\; \phi_3(t)=t^2(1+\tfrac54t+\ldots) \] 
Consequently monodromy for $H^{3}(Y_{t}), t\in\Delta^{*}$ has order $2$. 
\hfill $\Diamond$.

\section{Cohomological analysis}

To describe the cohomological relation between the singular 
fibre $Y_0$ and the generic fibre $Y_t$ of our family $f: \YS \to \Delta$, 
we use the nearby and vanishing cycle formalism from \cite{SGA7II}. 
There is a distinguished triangle in $D^b_{constr}(Y_0,\Q)$ that reads 
\[ \ldots \lra \Q_{Y_0} \lra R\Psi_f(\Q) \lra R\Phi_f(\Q) \stackrel{+1}{\lra} \ldots\]
and which leads to a long exact sequence in cohomology
\[\ldots \lra H^k(Y_0,\Q) \lra \H^k(R\Psi_f(\Q)) \lra \H^k(R\Phi_f(\Q)) \lra \ldots\] 
The cohomology group $H^k(Y_0,\Q)$ carries, after Deligne \cite{Deligne}, 
a natural mixed Hodge structure. The hypercohomology $\H^k(R\Psi_f(\Q))$ can be identified with the cohomology $H^k(Y_{\infty},\Q)$ of the nearby fibre and carries a canonical mixed Hodge  structure after Schmid \cite{Schmid} and Steenbrink \cite{Steenbrink}. 
The vanishing cohomology groups

$$\H^k:=\H^k(R\Phi_f(\Q))$$

can be given a mixed Hodge structure in a way compatible with this exact sequence.\\

{\bf Proposition:} One has 
\[ \H^k=0\;\;\textup{for}\;\;k \neq 3\]
and there is a short exact sequence
\[ 0 \lra H^1(R^2\Phi_f(\Q)) \lra \H^3 \lra H^0(R^3\Phi_f(\Q)) \lra 0\]
and identifications
\[ H^1(R^2\Phi_f(\Q))=H^1(E)(-1),\;\;\;H^0(R^3\Phi_f(\Q))=\oplus_{p\in P} \Q(-2).[p]\]
as MHS.\\

{\bf proof:} We use the hypercohomology spectral sequence $\H^p(R\Phi^q_f(\Q)) \implies \H^{p+q}$. As the singular locus is $L$, which is codimension $2$ in $Y_0$, we
have that $R^0\Phi_f(\Q)=R^1\Phi_f(\Q)=0$. At a general point $q \in L$ the 
threefold $Y_0$ has a transverse $A_1$-singularity, hence we have that the
stalk $R^2\Phi_f(\Q)_q $ is one-dimensional, whereas $R^3\Phi_f(\Q)_q=0$. 
At the pinch-points $p \in \Sigma$ one has $R^2\Phi_f(\Q)_p=0$ and 
$R^3\Phi_f(\Q)_p$ one-dimensional. So $R^3\Phi_f(\Q)$ is a sky-scraper sheaf
at the pinch points $\Sigma$, whereas $R^2\Phi_f(\Q)$ is a rank one-local 
system on $L\setminus \Sigma$, extended to zero.
From the local normal form of the pinch-point, one sees that the monodromy 
on $R^2\Phi_f(\Q)$ around the points $p \in \Sigma$ are multiplication by $(-1)$. Hence there is an exact sequence
\[ 0 \lra \Q_L \lra n_*(\Q_E) \lra R^2\Phi(\Q)(1) \lra 0,\] 
where $n:E \lra L$ is the elliptic curve over $L$, ramified at $\Sigma$.
From the long exact cohomology sequence we immediately obtain
\[ H^0(R^2\Phi_f(\Q))=0,\;\;H^1(R^2\Phi_f(\Q))=H^1(E,\Q)(-1),\;\; H^2(R^2\Phi_f(\Q))=0 .\]
The hypercohomology spectral sequences collapses and give the
above result. \hfill $\Diamond$\\

{\bf Corollary:}
The Hodge structures $H^k(Y_0,\Q)$ are pure of weight $k$. 
There are isomorphisms 
\[ H^k(Y_0,\Q) \approx H^k(Y_{\infty},\Q)\;\;\textup{for}\;k\neq 3, 4\]
and short exact sequences
\[ 0 \lra H^3(Y_0,\Q) \lra H^3(Y_{\infty},\Q) \lra H^1(E,\Q)(-1) \lra 0\]
\[ 0 \lra \oplus_{p \in \Sigma} \Q(-2)p \lra H^4(Y_{0},\Q) \lra H^4(Y_{\infty},\Q) \lra 0 .\]
Hence we have:
\[ h^2(Y_0)=h^2(Y_{\infty})=41,\;\; h^3(Y_0)=2,\;\;h^4(Y_0)=h^4(Y_{\infty})+4=45.\]
{\bf proof:} Note that the Hodge structures $H^k(Y_{\infty},\Q)$ are pure, 
as the monodromy is trivial. As only $\H^3 \neq 0$, we get from  
the long exact cohomology sequence isomorphisms and a five term exact sequence
\[ 0 \lra H^3(Y_0,\Q) \lra H^3(Y_{\infty},\Q) \lra \H^3 \lra H^4(Y_0,\Q) \lra H^4(Y_{\infty},\Q)\lra 0\]
From the fact that $H^3(Y_0,\Q) \neq H^3(Y_{\infty},\Q)$ and the fact that $H^3(Y_{\infty},\Q)$ is pure, we see the only possibility is that $H^3(Y_{\infty},\Q)$ 
surjects on the weight $3$ part $H^1(E,\Q)(-1)$ of $\H^3$ and then the kernel 
of  $H^4(Y_0,\Q) \lra H^4(Y_{\infty},\Q)$ is equal to the weight $4$ quotient of 
$\H^3$.
\hfill $\Diamond$.\\

\subsection{Cohomology of $Z_0$}
The space $Z_0$ is obtained from $Y_0$ by a single blow-up in the line
$L$. The preimage of $L$ is a conic-bundle $Q \lra L$. So
we get a diagram
\[
\begin{array}{ccc}
Q&\hookrightarrow&Z_0\\
\downarrow&&\downarrow \pi\\
L& \hookrightarrow&Y_0\\
\end{array}
\]

The following calculation provides an independent determination of
the cohomology of $Z_0$.\\

{\bf Proposition:} There are exact sequences of mixed Hodge structures
\[ 0 \lra H^2(Y_0,\Q) \lra H^2(Z_0,\Q) \lra H^0(R^2\pi_*\Q_{Z_0}) \lra 0.\]
\[  H^3(Y_0,\Q) \approx H^3(Z_0,\Q),\]
\[ 0 \lra H^4(Y_0,\Q) \lra H^4(Z_0,\Q) \lra H^2(R^2\pi_*\Q_{Z_0}) \lra 0 .\] 
Furthermore:
\[ H^0(R^2\pi_*\Q_{Z_0})=\Q(-1)^5,\;\;\; H^2(R^2\pi_*\Q_{Z_0})=\Q(-2),\]
hence one gets
\[ h^2(Z_0)=h^2(Y_0)+5=46,\;\;h^3(Z_0)=h^3(Y_0)=2,\;\;h^4(Z_0)=h^4(Y_0)+1=46 .\]

{\bf proof:}
In order to compare the cohomology of $Z_0$ and $Y_0$, we use the
Leray spectral sequence
\[ E_2^{p,q}=H^p(R^q\pi_*\Q_{Z_0}) \implies H^{p+q}(Z_0,\Q).\] 
The inverse image of $L$ in $Y$ is a conic bundle $Q \lra L$. 
The general fibre is a smooth conic; the fibres over $p \in \Sigma$ are line pairs. Hence
\[R^0\pi_*\Q_{Z_0}=\Q_{Y_0},\;\;\;R^k\pi_*\Q_{Z_0}=0,\;\;\;k \neq 0,2\]
From the study of the vanishing cohomology for the family $Q \lra L$ near 
the line-pairs, we obtain a split exact sequence
\[ 0 \lra \Q_{\Sigma} \lra R^2\pi_*\Q_Q \lra \Q_L \lra 0,\]
which indeed leads to 
\[ H^0(Q)=\Q,\;\;\;H^1(Q,\Q)=0,\;\;\;H^2(Q,\Q)=\Q(-1)^6,\;\;\;H^3(Q,\Q)=0,\;\;\;H^4(Q,\Q)=\Q(-2).\]
The differential
\[ d_3:H^0(R^2\pi_*\Q_{Z_0}) \to H^3(R^0\pi_*\Q_{Z_0})=H^3(Y_0,\Q)\]
has to be zero because of weights, so we obtain from the spectral
sequence short exact sequences as stated above.
\hfill $\Diamond$.

\subsection{The semi-stable reduction}
Denote by
$g: \ZS \to \YS$ the blow-up of the smooth space $\YS$ in the line $L$.
We denote the exceptional divisor of this blow-up by $P$; it is a 
$\P^2$-bundle over $L$. As the multiplicity of $Y_0$ alone $L$ is two,
the divisor of the composed function 
\[ h:=f\circ g: \ZS \to \Delta\]
is 
\[ Z_0 + 2 P,\]
where the strict transform $Z_0$ of $Y_0$ is blow-up of $Y_0$ in $L$,
hence smooth. The intersection of these two components is the surface
\[ Q := Z_0 \cap P .\]
The map $Q \to L$, obtained as restriction of $g$, gives $Q$ the structure
of a conic bundle over $L$; above the four pinch-points the conics
degenerate into a line pair.
Now we take the pull-back of $h:\ZS \to \Delta$ by the squaring map 
$s:t \mapsto t^2$, and denote its normalisation by $\widetilde{\ZS}$:
\[\widetilde{\ZS}:=\widetilde{\Delta \times_{\Delta} \ZS} .\]
We let $n:\widetilde{\ZS}  \to \ZS$ the natural map, so we have the
diagram

\[
  \begin{diagram}
    \dgARROWLENGTH=8mm
    \node{\widetilde{P}}\arrow[2]{s}\node{+}\node{\widetilde{Z_0}}\arrow[2]{s}\node{\subset}
  \node{\widetilde{\ZS}}\arrow[2]{s,l}{\tilde g} \arrow[2]{e,t}{n}
  \node[2]{\ZS}\arrow[2]{s,r}{g}\node{\supset}\node{Z_0}\arrow[2]{s}\node{+}\node{P}\arrow[2]{s}\\[2]
  \node{\widetilde{L}}\node{\subset}\node{\widetilde{\YS_0}}
  \node{\subset} \node{\widetilde{\YS}}\arrow[2]{s}\arrow[2]{e}
  \node[2]{\YS}\arrow[2]{s}\node{\supset}\node{Y_0}\node{\supset} \node{L} \\[2]
  \node[5]{\Delta}\arrow[2]{e,t}{s} \node[2]{\Delta}
  \end{diagram}
\]

{\bf Proposition:} The space $\widetilde{\ZS}$ is smooth.
The divisor 
\[\widetilde{h}^{-1}(0) = \widetilde{Z_0} + \widetilde{P}\]
is reduced and normal crossing. 
The map $n$ induces an isomorphism
\[ \widetilde{Z_0} \to Z_0\]
and a $2$-to-$1$ covering $\widetilde{P} \to P$ ramified precisely along $Q \subset P$.\\
The cohomology groups $H^i(\widetilde{P},\Q)$ are given by
\[\def\arraystretch{1.2}
\begin{array}{|c|c|c|c|c|c|c|}
\hline
H^0&H^1&H^2&H^3&H^4&H^5&H^6\\[1mm]
\hline
\Q&0&\Q(-1)^2&H^1(E)(-1)&\Q(-2)^2&0&\Q(-3)\\[1mm]
\hline
\end{array}
\]

{\bf proof:} This follows from a direct local calculation.
Around any point of $Q=Z_0 \cap P$ the divisor $Z_0+2P$ is
given by an equation of the form $xy^2=0$. The $2$-fold cover
then has equation $xy^2+z^2=0$, which has a smooth normalisation.   
Clearly, the map $\widetilde{P} \to P$ is a two fold cover,  ramified 
precisely along the conic bundle $Q$. In other words, the composition
$\rho: \widetilde{P} \to P \to L$ represents this threefold as a quadric
bundle, with four fibres with an isolated singular point over the points 
$\Sigma$. We can determine the cohomology of $\widetilde{P}$ using the
Leray spectral sequence of the map $\rho:\widetilde{P} \to L$. 
We find
\[ R^0\rho_*(\Q_{\widetilde{P}})=\Q_L,\;\;R^1\rho_*(\Q_{\widetilde{P}})=0=R^3\rho_*(\Q_{\widetilde{P}}),\;\;R^4\rho_*(\Q_{\widetilde{P}})=\Q(-2) .\]
The sheaf $R^2\rho_*(\Q_{\widetilde{P}})$ is more interesting. As $H^2$ of
of a quadric is generated by its two rulings, which get interchanged upon
surrounding a point of $\Sigma$, and coalesce over $\Sigma$, we have
\[ R^2\rho_*(\Q_{\widetilde{P}})=\pi_*\Q_E(-1),\]
where $\pi:E \lra L$ is the elliptic curve, two-fold covering $L$ and
ramifying over $\Sigma$. The Leray-spectral sequence degenerates and we
can read off directly the cohomology groups, as Hodge structures. 
\hfill $\Diamond$.\\

The monodromy weight spectral sequence converges to the cohomology of
$H^k(\widetilde{Z}_{\infty},\Q)=H^k(Y_{\infty},\Q)$ and is determined from 
the intersections of the irreducible  components of the semi-stable fibre,
see e.g. \cite{Steenbrink}, \cite{SteenbrinkPeters}. 
In our case there are only two components and a single intersection, so the 
$E_1^{p,q}$-page is very simple and looks like
\[
\begin{array}{ccccc}
H^4(Q)(-1) &\lra&H^6({Z}_0)\oplus H^6(\widetilde{P})&\lra&0\\
H^3(Q)(-1) &\lra&H^5({Z}_0)\oplus H^5(\widetilde{P})&\lra&0\\
H^2(Q)(-1) &\lra&H^4({Z}_0)\oplus H^4(\widetilde{P})&\lra&H^4(Q)\\
H^1(Q)(-1) &\lra&H^3({Z}_0)\oplus H^3(\widetilde{P})&\lra&H^3(Q)\\
H^0(Q)(-1) &\lra&H^2({Z}_0)\oplus H^2(\widetilde{P})&\lra&H^2(Q)\\
0          &\lra&H^1({Z}_0)\oplus H^1(\widetilde{P})&\lra&H^1(Q)\\
0          &\lra&H^0({Z}_0)\oplus H^0(\widetilde{P})&\lra&H^0(Q)\\
\end{array}
\]
As we have determined all groups appearing, the cohomology-diagram 
has the following form:
\[
\begin{array}{ccccc}
\Q(-3)  &\lra &\Q(-3) \oplus \Q(-3)       &\lra&0\\
0       &\lra & 0     \oplus 0            &\lra&0\\
\Q(-2)^6&\lra & \Q(-2)^{46} \oplus \Q(-2)^2&\lra&\Q(-2)\\
0       &\lra &  H^3(Z_0,\Q) \oplus H^1(E,\Q)(-1) &\lra&0\\
\Q(-1)  &\lra&\Q(-1)^{46} \oplus \Q(-1)^2&\lra&\Q(-1)^6\\
0       &\lra&0 \oplus  0               &\lra&0\\
0       &\lra&\Q \oplus \Q              &\lra&\Q\\
\end{array}
\]

As we know that the cohomology groups $H^k(Y_{\infty},\Q)$ of the limit are
pure Hodge structures, the maps at the left are injective, those on the right 
surjective, and the cohomology of $H^k(Y_{\infty},\Q)$ comes out the right way:
\[H^3_{\lim}(Y_0,\Q)= H^3(Z_0,\Q) \oplus H^1(E,\Q)(-1).\]

\section{Outlook}

The family $\mathcal X\lra\Delta$ extends naturally to a
projective family over $\mathbb P^{1}$ with four singular fibers
at $\{0,1,2,\infty\}$. The fibers at $1$ and $2$ are double octic
arrangements No. $3$ and No. $19$ respectively, whereas at $\infty$ 
we get degenerate configuration
\[u^{2}=xy(x+y)zv^{3}(y+z+v).\]
The map
\[(x,y,z,v,u)\longmapsto
  (-x,x+y,z,\tfrac1{t-1}v,u)\]
defines an isomorphism between $\XS_{\frac t{t-1}}$ and the quadratic
twist of $\XS_{t}$ by $1-t$, consequently quadratic base-changes of the
family $\XS$ ramified at $1$ and $\infty$ are isomorphic.

There are $63$ one-parameter families of double octics listed in \cite{Meyer},
which lead to $63$ families of Calabi-Yau threefolds with $h^{12}=1$. In these
families there are five more examples with a similar behaviour. Below one
list all six cases, which come as three pairs:\\

\noindent
$$
\begin{array}{|r|l|c|r|}
\hline
No.& equation&t_{0}&\XS_{t_{0}}\\
\hline
\hline
153&xyzt \left( x+y+z \right)  \left( y+z+t
 \right)\times &-2&93\\
& \times   \left( Ax-By+At
\right) \left( Ax-By+Az+At \right)&&\\
\hline
 197&xyzt(x-y-z+t)(Ax+By+Bz)\times &-\frac12& 93\\
 &\times (By+Bz+At)(Ax+Bz+At)&&\\
\hline
\hline
96&xyzt(x+y)(x+y-z+t)\times&-2&32\\
 &\times(Ax-By+Bz+At)(Ay+Bz+At)&&\\
\hline
 100&xyzt(x+y-z+t)(Ax+Ay+Bz)\times&-\frac12&69\\
 &\times(Ay+Bz+At)(By-Bz-At)&&\\
\hline
\hline
155&xyzt(Ax+By+Az)(Ax+(A+B)y-Bz+At)\times &
 \frac{-1\pm\sqrt{-3}}2&A\\
 &\times (Ax-Bz-Bt)(Ax+By+Az+At)&&\\
\hline
200&xyzt(x+y+z+t)(Ax+Ay-Bz-Bt)\times &\frac{-1\pm\sqrt{-3}}2&A\\
 &\times (Ay-Bz+At)(Ax-By-Bt)&&\\
\hline
 \end{array}
$$
\vskip 10pt
In the last column we have indicated the configuration number of the
corresponding double octic from \cite{Meyer}. The symbol $A$ indicates
a specific rigid Calabi-Yau manifold defined over $\Q(\sqrt{-3})$.
The families No.  $96$ and No. $100$ are in fact birational, 
as are No. $153$ and No. $197$. Families No. $155$ and No. $200$
have equal Hodge numbers and share the same Picard-Fuchs operators, 
but no birational map between them is known to us. 

The degeneration of two fourfold points of type $p_{4}^{1}$ that 
collide and produce a $p_{5}^{1}$ that was analysed in this paper
for No. $153$ also occurs in No. $100$ and No. $155$. As a 
consequence we get again the central fiber singular along a double 
line with four pinch points. The only difference is that in the 
case of family No. $155$ the $j$-invariant of four pinch-points equals $0$.
The degenerations that occur in the other three cases No. $96$, No. $197$ 
and No. $200$ are of a different kind: three double lines come together 
to form a triple line.  This line is a double line of the singular element 
of the central fiber with four pinch point: one fivefold point and three 
fourfold points on this line. The $j$-invariant of this four points is again 
$1728$ in first two cases and $0$ in the last case. 

The local exponents of the Picard-Fuchs operators in the first four
families are all  equal to $(0,1/2,3/2,2)$, which after quadratic base 
change become $(0,1,3,4)$, while in the case 
of the last two families No. $155, 200$ they are $(0,1/2,5/2,3)$, which
after a quadratic base change  become $(0,1,5,6)$. 
It is surprising and beautiful to see the order of the automorphism group 
of the associated elliptic curve appear in the local exponents of the 
degeneration.

The first four families No. $96, 100, 153$ and $197$ are also birational 
to Kummer fibrations of rational elliptic surfaces. The degeneration of 
the corresponding fiber products results from the collisions of fibers

\[
    (I_{2}\times I_{0}^{*})+(I_{0}\times I_{0})
\lra    (I_{2}\times I_{0}^{*})
  \qquad\text{ or }\qquad
    (I_{2}\times I_{0}^{*})+(I_{2}\times I_{0})
\lra    (I_{4}\times I_{0}^{*})
\]
The singularities in the central fiber correspond to two copies of
the singular fiber of type $I_{0}^{*}$. It may very well be possible
to analyse the degeneration cohomologically from this description.\\

We believe that our degeneration also has an arithmetical version
that may be of interest. Recently, in \cite{LiedtkeMatsumoto} and \cite{CLL} 
a version N\'eron-Ogg-Shafarevich criterion for a family of K3 surfaces
was formulated, which enables to detect good reduction of a K3
surface over the fraction field $K$ of a henselian local ring $R$ 
with residue field $k$ of characteristic $p >0$ by having the
Galois representation
\[ G_K \to Aut(H_{et}^2(X,\Q_{\ell}))\]
unramified. Based on our example, we are inclined to think that no 
similar criterion can exist for Calabi-Yau threefolds. If we replace 
$t$ by a sufficiently large prime $p$ in the formula describing our 
double octic, and doing the corresponding modifications, we end up 
with with a Calabi-Yau variety $Y$ over the $p$-adic field $K=\Q_p$ 
for which we have the suspicion that the Galois representation
\[ G_K \to Aut(H_{et}^3(Y,\Q_{\ell}))\]
is unramified for $\ell \neq p$ and crystalline for $\ell=p$, but for
which no good (terminal) reduction is in sight.\\

{\bf Acknowledgement:} We thank Radu Laza for showing interest in this
example. The first named author was partially supported by the National 
Science Center grant no. 2014/13/B/ST1/00133. This work was partially 
supported by the grant 346300 for IMPAN from the Simons Foundation and the 
matching 2015-2019 Polish MNiSW fund.

\end{document}